\documentclass[11pt]{amsart}
\usepackage{amsmath, amsthm, amsfonts, amssymb} 
\usepackage{url}

\newtheorem{theo}{Theorem}[section]
\newtheorem*{theorem}{Theorem}

\newtheorem{cor}[theo]{Corollary}
\newtheorem{lem}[theo]{Lemma}

\theoremstyle{definition}

\newtheorem{nota}[theo]{Notation}
\newtheorem{df}[theo]{Definition}

\newcommand{\R}{\mathbf{R}}
\newcommand{\C}{\mathbf{C}}

\newcommand{\F}{\mathbf{F}}

\newcommand{\N}{\mathbf{N}}

\newcommand{\Sp}{\operatorname{Sp}}

\newcommand{\Ad}{\operatorname{Ad}}
\newcommand{\id}{\operatorname{id}}
\newcommand{\Tr}{\operatorname{Tr}}

\newcommand{\AC}{\operatorname{AC}}
\newcommand{\AB}{\operatorname{AB}}
\newcommand{\dom}{\operatorname{dom}}
\newcommand{\support}{\operatorname{support}}

\begin{document}

\title[Free Araki-Woods Factors, Bicentralizer Problem]{Free Araki-Woods Factors \\ and Connes' Bicentralizer Problem}

\begin{abstract}
We show that for any type ${\rm III_1}$ free Araki-Woods factor $\mathcal{M} = \Gamma(H_\R, U_t)''$, the bicentralizer of the free quasi-free state $\varphi_U$ is trivial. Using Haagerup's Theorem, it follows that there always exists a faithful normal state $\psi$ on $\mathcal{M}$ such that $(\mathcal{M}^\psi)' \cap \mathcal{M} = \C$.
\end{abstract}

\author{Cyril Houdayer }

\address{CNRS-ENS Lyon \\
UMPA UMR 5669 \\
69364 Lyon cedex 7 \\
France}

\email{cyril.houdayer@umpa.ens-lyon.fr}

\subjclass[2000]{46L10; 46L54}

\keywords{Free Araki-Woods factors; Connes' bicentralizer problem}

\maketitle

\section{Introduction}

Let $\mathcal{M}$ be a separable type ${\rm III_1}$ factor and let $\varphi$ be a faithful normal state on $\mathcal{M}$. For any $x, y \in \mathcal{M}$, set $[x, y] = xy - yx$ and $[x, \varphi] = x \varphi - \varphi x$. The {\it asymptotic centralizer} of $\varphi$ is defined by
\begin{equation*}
\AC(\varphi) := \{(x_n) \in \ell^\infty(\N, \mathcal{M}) : \|[x_n, \varphi]\| \to 0\}.
\end{equation*}
Note that $\AC(\varphi)$ is a unital C$^*$-subalgebra of $\ell^\infty(\N, \mathcal{M})$. The {\it bicentralizer} of $\varphi$ is defined by
\begin{equation*}
\AB(\varphi) := \{a \in \mathcal{M} : [a, x_n] \to 0 \; \mbox{ultrastrongly}, \forall (x_n) \in \AC(\varphi)\}.
\end{equation*}
It is well-known that $\AB(\varphi)$ is a von Neumann subalgebra of $\mathcal{M}$, globally invariant under the modular group $(\sigma_t^\varphi)$. Moreover, $\AB(\varphi) \subset (\mathcal{M}^\varphi)' \cap \mathcal{M}$. If $\AB(\varphi) = \C$, it follows from Connes-St\o rmer Transitivity Theorem (\cite{connesstormer}), that $\AB(\psi) = \C$ for any faithful normal state $\psi$ on $\mathcal{M}$. We shall say in this case that $\mathcal{M}$ has {\it trivial bicentralizer}. Connes conjectured that {\it any} separable type ${\rm III_1}$ factor should have trivial bicentralizer. If there exists a faithful normal state $\varphi$ on $\mathcal{M}$ such that $(\mathcal{M}^\varphi)' \cap \mathcal{M} = \C$, then $\mathcal{M}$ has trivial bicentralizer. Haagerup proved in \cite{haagerup84} that the converse holds true. Haagerup's Theorem lead to the uniqueness of the amenable type ${\rm III_1}$ factor (see \cite{connes85}). The following type ${\rm III_1}$ factors are known to have trivial bicentralizer:

\begin{enumerate}
\item The unique amenable ${\rm III_1}$ factor (Haagerup, \cite{haagerup84}).
\item Full factors that have almost periodic states (Connes, \cite{connes74}).
\item Free products $(\mathcal{M}_1, \varphi_1) \ast (\mathcal{M}_2, \varphi_2)$ such that the centralizers $\mathcal{M}_i^{\varphi_i}$ have enough unitaries (Barnett, \cite{barnett95}).
\end{enumerate}

In this paper, we show that the bicentralizer is trivial for a large class of type ${\rm III_1}$ factors, namely the {\it free Araki-Woods factors} of Shlyakhtenko (\cite{shlya97}). We briefly recall the construction here, see Section \ref{preliminaries} for more details. To each real separable Hilbert space $H_\R$ together with an orthogonal representation $(U_t)$ of $\R$ on $H_\R$, one can associate a von Neumann algebra denoted by $\Gamma(H_\R, U_t)''$, called the free Araki-Woods von Neumann algebra. This is the free analog of the factors coming from the CAR relations. The von Neumann algebra $\Gamma(H_\R, U_t)''$ comes equipped with a unique {\it free quasi-free state} denoted by $\varphi_U$, which is always normal and faithful on $\Gamma(H_\R, U_t)''$. If $\dim H_\R \geq 2$, then $\Gamma(H_\R, U_t)''$ is a full factor. It is of type ${\rm III_1}$ when $(U_t)$ is non-periodic and non-trivial. If the representation $(U_t)$ is almost periodic, then $\varphi_U$ is an almost periodic state and it follows from \cite{shlya97} that the relative commutant of the centralizer of the free quasi-free state is trivial, i.e. if $\mathcal{M} := \Gamma(H_\R, U_t)''$, then $(\mathcal{M}^{\varphi_U})' \cap \mathcal{M} = \C$. In the almost periodic case, results of \cite{dykema96} yield $\mathcal{M}^{\varphi_U} \simeq L(\F_\infty)$.

When the representation $(U_t)$ has no eigenvectors (e.g. $U_t = \lambda_t$,  the left regular representation of $\R$ on $L^2(\R, \R)$), then the centralizer $\mathcal{M}^{\varphi_U}$
 is trivial. It was unknown in general whether or not $\Gamma(H_\R, U_t)''$ had trivial bicentralizer. Even though the centralizer of the free quasi-free state $\varphi_U$ may be trivial, we will show that the bicentralizer of $\varphi_U$ is always trivial. The main result of this paper is the following:

\begin{theorem}%\label{bicentralizer}
Let $\mathcal{M} := \Gamma(H_\R, U_t)''$ be a free Araki-Woods factor of type ${\rm III_1}$. Denote by $\varphi_U$ the free quasi-free state. Then $\AB(\varphi_U) = \C$. Consequently, there always exists a faithful normal state $\psi$ on $\mathcal{M}$ such that $(\mathcal{M}^\psi)' \cap \mathcal{M} = \C$. 
\end{theorem}

{\bf Acknowledgements.} Part of this work was done while the author was visiting the University of Tokyo. He gratefully thanks Professors Y. Kawahigashi and N. Ozawa for their kind invitation.

\section{Preliminaries}\label{preliminaries}

\subsection{Preliminaries on Spectral Analysis}

We shall need a few definitions and results from the spectral theory of abelian automorphism groups. Let $(\alpha_t)$ be an ultraweakly continuous one-parameter automorphism group on a von Neumann algebra $\mathcal{M}$. For $f \in L^1(\R)$ and $x \in \mathcal{M}$, set
\begin{equation*}
\alpha_f(x) = \int_{- \infty}^{+ \infty} f(t) \alpha_t(x) \,dt.
\end{equation*}
The $\alpha$-spectrum $\Sp_{\alpha}(x)$ of $x \in \mathcal{M}$ is defined as the set of characters $\gamma \in \widehat{\R}$ for which $\widehat{f}(\gamma) = 0$, for all $f \in L^1(\R)$ satisfying $\alpha_f(x) = 0$. We shall identify $\widehat{\R}$ with $\R$ in the usual way, such that
\begin{equation*}
\widehat{f}(\gamma) = \int_{- \infty}^{+ \infty} e^{i \gamma t}f(t) \,dt, \forall \gamma \in \R, \forall f \in L^1(\R).
\end{equation*}
For $z \in \C$, denote by $\Im(z)$ its imaginary part.

\begin{lem}[\cite{haagerup84}]\label{spectrum1}
Let $\mathcal{M}$ and $(\alpha_t)$ be as above. Let $x \in \mathcal{M}$ and $\delta > 0$. If the function $t \mapsto \alpha_t(x)$ can be extended to an entire (analytic) $\mathcal{M}$-valued function such that
\begin{equation*}
\|\alpha_z(x)\| \leq C e^{\delta |\Im(z)|}, \forall z \in \C,
\end{equation*}
for some constant $C > 0$, then $\Sp_{\alpha}(x) \subset [-\delta, \delta]$.
\end{lem}

Let $\varphi$ be a f.n. state on a von Neumann algebra $\mathcal{M}$. Denote by $(\sigma^\varphi_t)$ the modular group on $\mathcal{M}$ of the state $\varphi$. Denote by $L^2(\mathcal{M}, \varphi)$ the $L^2$-space associated with $\varphi$ and by $\xi_\varphi$ the canonical cyclic separating vector. 
We shall write $\|x\|_\varphi = \varphi(x^*x)^{1/2}$, for any $x \in \mathcal{M}$. On bounded subsets of $\mathcal{M}$, the topology given by the norm $\|\cdot\|_\varphi$ coincides with the strong-operator topology. Recall that $S^0_\varphi : x \xi_\varphi \mapsto x^*\xi_\varphi$ is a closable (densely defined) operator on $L^2(\mathcal{M}, \varphi)$. Denote by $S_\varphi$ its closure and write $S_\varphi = J_\varphi \Delta^{1/2}_\varphi$ for its polar decomposition. Note that $L^2(\mathcal{M}, \varphi)$ is naturally endowed with an $\mathcal{M}$-$\mathcal{M}$ bimodule structure defined as follows:
\begin{eqnarray*}
x \cdot \xi & := & x\xi \\
\xi \cdot x & := & J_\varphi x^* J_\varphi \xi, \forall x \in \mathcal{M}, \forall \xi \in L^2(\mathcal{M}, \varphi).
\end{eqnarray*}
We shall simply denote $x \cdot \xi$ and $\xi \cdot x$ by $x \xi$ and $\xi x$. The next lemma is well-known, but we give a proof for the reader's convenience.

\begin{lem}\label{spectrum2}
Let $\mathcal{M}$ and $\varphi$ be as above. Let $x \in \mathcal{M}$ and $0 < \delta < 1$. Assume that $\Sp_{\sigma^\varphi}(x) \subset [-\delta, \delta]$. Then $\|x \xi_\varphi - \xi_\varphi x\| \leq \delta \|x\|_\varphi$. 
\end{lem}

\begin{proof}
Let $x \in \mathcal{M}$ and $0 < \delta < 1$ such that $\Sp_{\sigma^\varphi}(x) \subset [-\delta, \delta]$. Let $f \in L^1(\R)$ for which the Fourier tranform $\widehat{f}$ vanishes on $[-\delta, \delta]$. Since $\Sp_{\sigma^\varphi}(\sigma^\varphi_f(x)) \subset \Sp_{\sigma^\varphi}(x) \cap \support(\widehat{f}) = \varnothing$ (see \cite{connes73}), it follows that $\sigma^\varphi_f(x) = 0$. We have
\begin{eqnarray*}
\widehat{f}(\log \Delta_\varphi) x\xi_\varphi & = & \int_{-\infty}^{+\infty} f(t) \Delta_\varphi^{it} x\xi_\varphi \,dt \\
& = & \int_{-\infty}^{+\infty} f(t) \sigma^\varphi_t (x)\xi_\varphi \,dt \\
& = & \sigma^\varphi_f(x) \xi_\varphi \\
& = & 0.
\end{eqnarray*}
Thus, by approximating ${\bf 1}_{\R \backslash [-\delta, \delta]}$ by such functions $\widehat{f}$, we get 
\begin{equation*}
{\bf 1}_{\R \backslash [-\delta, \delta]}(\log\Delta_\varphi)x\xi_\varphi = 0,
\end{equation*}
i.e. $x \xi_\varphi$ is in the spectral subspace of $\log \Delta_\varphi$ corresponding to the interval $[-\delta, \delta]$. Notice that
\begin{equation*}
\xi_\varphi x = J_\varphi x^* J_\varphi \xi_\varphi = J_\varphi x^* \xi_\varphi = J_\varphi S_\varphi x \xi_\varphi = \Delta^{1/2}_\varphi x \xi_\varphi. 
\end{equation*}
Clearly, $\sup \{|e^{t/2} - 1| : t \in [-\delta, \delta]\} = e^{\delta/2} - 1$. Moreover, the operator $(1 - \Delta_\varphi^{1/2}) {\bf 1}_{[-\delta, \delta]}(\log\Delta_\varphi)$ is bounded and precisely 
\begin{equation*}
\|(1 - \Delta_\varphi^{1/2}) {\bf 1}_{[-\delta, \delta]}(\log\Delta_\varphi)\| \leq e^{\delta/2} - 1 \leq \delta,
\end{equation*}
since $0 < \delta < 1$. Thus, we get
\begin{eqnarray*}
\|x \xi_\varphi - \xi_\varphi x\|  & = & \|(1 - \Delta_\varphi^{1/2})x \xi_\varphi\| \\
& = & \|(1 - \Delta_\varphi^{1/2}){\bf 1}_{[-\delta, \delta]}(\log\Delta_\varphi)x \xi_\varphi\| \\
& \leq & \|(1 - \Delta_\varphi^{1/2}) {\bf 1}_{[-\delta, \delta]}(\log\Delta_\varphi)\| \|x \xi_\varphi\| \\
& \leq & \delta \|x\|_\varphi.
\end{eqnarray*}
\end{proof}

\begin{lem}[\cite{haagerup84}]\label{spectrum3}
Let $\mathcal{M}$ and $\varphi$ be as above. Let $(x_n) \in \ell^\infty(\N, \mathcal{M})$. Then
\begin{equation*}
\lim_n \|x_n \xi_\varphi - \xi_\varphi x_n\| = 0 \Longleftrightarrow \lim_n \|x_n \varphi - \varphi x_n\| = 0.
\end{equation*}
\end{lem}

\subsection{Preliminaries on Shlyakhtenko's Free Araki-Woods Factors}

Recall now the construction of the free Araki-Woods factors due to Shlyakhtenko (\cite{shlya97}). Let $H_{\R}$ be a real separable Hilbert space and let $(U_t)$ be an orthogonal representation of $\R$ on $H_{\R}$. Let $H = H_{\R} \otimes_{\R} \C$ be the complexified Hilbert space. Let $J$ be the canonical anti-unitary involution on $H$ defined by:
\begin{equation*}
J(\xi + i \eta) = \xi - i \eta, \forall \xi, \eta \in H_\R.
\end{equation*}
 If $A$ is the infinitesimal generator of $(U_t)$ on $H$, we recall that $j : H_{\R} \to H$ defined by $j(\zeta) = (\frac{2}{A^{-1} + 1})^{1/2}\zeta$ is an isometric embedding of $H_{\R}$ into $H$. Moreover $JAJ = A^{-1}$ and $J A^{it} = A^{it} J$, for every $t \in \R$. Let $K_{\R} = j(H_{\R})$. It is easy to see that $K_\R \cap i K_\R = \{0\}$ and $K_\R + i K_\R$ is dense in $H$. Write $T = J A^{-1/2}$. Then $T$ is an anti-linear closed invertible operator on $H$ satisfying $T = T^{-1}$. Such an operator is called an {\it involution} on $H$. Moreover, $K_\R = \{ \xi \in \dom(T) : T\xi = \xi \}$.

We introduce the \emph{full Fock space} of $H$:
\begin{equation*}
\mathcal{F}(H) =\C\Omega \oplus \bigoplus_{n = 1}^{\infty} H^{\otimes n}.
\end{equation*}
The unit vector $\Omega$ is called the \emph{vacuum vector}. For any $\xi \in H$, define the {\it left creation} operator
\begin{equation*}
\ell(\xi) : \mathcal{F}(H) \to \mathcal{F}(H) : \left\{ 
{\begin{array}{l} \ell(\xi)\Omega = \xi, \\ 
\ell(\xi)(\xi_1 \otimes \cdots \otimes \xi_n) = \xi \otimes \xi_1 \otimes \cdots \otimes \xi_n.
\end{array}} \right.
\end{equation*}
We have $\|\ell(\xi)\| = \|\xi\|$ and $\ell(\xi)$ is an isometry if $\|\xi\| = 1$. For any $\xi \in H$, we denote by $s(\xi)$ the real part of $\ell(\xi)$ given by
\begin{equation*}
s(\xi) = \frac{\ell(\xi) + \ell(\xi)^*}{2}.
\end{equation*}
The crucial result of Voiculescu \cite{voiculescu92} claims that the distribution of the operator $s(\xi)$ with respect to the vacuum vector state $\varphi(x) = \langle x\Omega, \Omega\rangle$ is the semicircular law of Wigner supported on the interval $[-\|\xi\|, \|\xi\|]$. 

\begin{df}[Shlyakhtenko, \cite{shlya97}]
Let $(U_t)$ be an orthogonal representation of $\R$ on the real Hilbert space $H_{\R}$. The \emph{free Araki-Woods} von Neumann algebra associated with $(H_\R, U_t)$, denoted by $\Gamma(H_{\R}, U_t)''$, is defined by
\begin{equation*}
\Gamma(H_{\R}, U_t)'' := \{s(\xi) : \xi \in K_{\R}\}''.
\end{equation*}
\end{df}

The vector state $\varphi_{U}(x) = \langle x\Omega, \Omega\rangle$ is called the {\it free quasi-free state} and is faithful on $\Gamma(H_\R, U_t)''$. Let $\xi, \eta \in K_\R$ and write $\zeta = \xi + i \eta$. We have
\begin{equation*}
2 s(\xi) + 2 i s(\eta) = \ell(\zeta) + \ell(T \zeta)^*.
\end{equation*}
Thus, $\Gamma(H_\R, U_t)''$ is generated as a von Neumann algebra by the operators of the form $\ell(\zeta) + \ell(T \zeta)^*$ where $\zeta \in \dom(T)$. Note that the modular group $(\sigma_t^{\varphi_U})$ of the free quasi-free state $\varphi_U$ is given by $\sigma^{\varphi_U}_{- t} = \Ad(\mathcal{F}(U_t))$, where $\mathcal{F}(U_t) = \id \oplus \bigoplus_{n \geq 1} U_t^{\otimes n}$. In particular, it satisfies
\begin{equation*}
\sigma_{-t}^{\varphi_U}\left( \ell(\zeta) + \ell(T \zeta)^* \right)  =  \ell(U_t \zeta) + \ell(T U_t \zeta)^*, \forall \zeta \in \dom(T), \forall t \in \R. 
\end{equation*}

The free Araki-Woods factors provided many new examples of full factors of type {\rm III} \cite{{barnett95}, {connes73}, {shlya2004}}. We can summarize the general properties of the free Araki-Woods factors in the following theorem (see also \cite{vaes2004}):

\begin{theo}[Shlyakhtenko, \cite{{shlya2004}, {shlya99}, {shlya98}, {shlya97}}]
Let $(U_t)$ be an orthogonal representation of $\R$ on the real Hilbert space $H_{\R}$ with $\dim H_{\R} \geq 2$. Denote by $\mathcal{M} := \Gamma(H_{\R}, U_t)''$.
\begin{enumerate}
\item $\mathcal{M}$ is a full factor and Connes' invariant $\tau(\mathcal{M})$ is the weakest topology on $\R$ that makes the map $t \mapsto U_t$ strongly continuous.
\item $\mathcal{M}$ is of type ${\rm II_1}$ iff $U_t = \id$, for every $t \in \R$.
\item $\mathcal{M}$ is of type ${\rm III_{\lambda}}$ $(0 < \lambda < 1)$ iff $(U_t)$ is periodic of period $\frac{2\pi}{|\log \lambda|}$.
\item $\mathcal{M}$ is of type ${\rm III_1}$ in the other cases.
\item The factor $\mathcal{M}$ has almost periodic states iff $(U_t)$ is almost periodic.
\end{enumerate}
\end{theo}
Let $H_{\R} = \R^2$ and $0 < \lambda < 1$. Let
\begin{equation}\label{periodic}
U_t^\lambda = \begin{pmatrix}
\cos(t\log \lambda) & - \sin(t\log \lambda) \\
\sin(t\log \lambda) & \cos(t\log \lambda)
\end{pmatrix}.
\end{equation}

\begin{nota}[\cite{shlya97}]\label{Tlambda}
\emph{Denote by $(T_{\lambda}, \varphi_{\lambda}) := (\Gamma(H_{\R}, U_t)'', \varphi_U)$ where $H_{\R} = \R^2$ and $(U_t)$ is given by Equation $(\ref{periodic})$.}
\end{nota}
Using a powerful tool called the \emph{matricial model}, Shlyakhtenko was able to prove the following isomorphism
\begin{equation*}
(T_{\lambda}, \varphi_{\lambda}) \cong (\mathbf{B}(\ell^2(\N)), \psi_{\lambda}) \ast (L^{\infty}[-1, 1], \mu),
\end{equation*}
where $\psi_{\lambda}(e_{ij}) = \delta_{ij}\lambda^j(1 - \lambda)$, $i, j \in \N$, and $\mu$ is a non-atomic measure on $[-1, 1]$. The notation $\cong$ means a state-preserving isomorphism. He also proved that $(T_{\lambda}, \varphi_{\lambda})$ has the {\it free absorption property}, namely,
\begin{equation*}
(T_\lambda, \varphi_\lambda) \ast L(\mathbf{F}_\infty) \cong (T_\lambda, \varphi_\lambda).
\end{equation*}

\section{The Main Result}

\subsection{Technical Lemmas}

As we said before, the centralizer of the free quasi-free state may be trivial; this is the case for instance when the orthogonal representation $(U_t)$ on $H_\R$ has no eigenvectors. Nevertheless, the following lemma shows that for any free Araki-Woods von Neumann algebra, there exists a non-trivial sequence of unitaries $(u_n)$ in the asymptotic centralizer of the free quasi-free state $\varphi_U$.

\begin{lem}[Vaes, \cite{vaes2004}]\label{asymptotic}
Let $\mathcal{M} := \Gamma(H_\R, U_t)''$ be a free Araki-Woods von Neumann algebra. Denote by $\varphi$ the free quasi-free state and by $(\sigma_t)$ the modular group of the state $\varphi$. Then there exists a sequence of unitaries $(u_n)$ in $\mathcal{M}$ entire (analytic) w.r.t. $(\sigma_t)$ such that 
\begin{enumerate}
\item $\|\sigma_z(u_n) - u_n\| \to 0$ uniformly on compact sets of $\C$; 
\item $\varphi(u_n) \to 0$;
\item $(u_n) \in \AC(\varphi)$.
\end{enumerate}
\end{lem}

\begin{proof}
This lemma, with the exception of item $(3)$, is Vaes' result (see Lemma 4.3 in \cite{vaes2004}). Item $(3)$ was not observed by Vaes but is immediate from the construction using Lemmas \ref{spectrum1}, \ref{spectrum2}, \ref{spectrum3}.
\end{proof}

The following lemma is a generalization of Barnett's lemma (see \cite{barnett95}) which was itself a generalization of Murray \& von Neumann's $14 \varepsilon$ lemma. 

\begin{lem}[Vaes, \cite{vaes2004}]\label{14eps}
For $i = 1, 2$, let $(\mathcal{M}_i, \varphi_i)$ be a von Neumann algebra endowed with a f.n. state. Denote by $(\mathcal{M}, \varphi) = (\mathcal{M}_1, \varphi_1) \ast (\mathcal{M}_2, \varphi_2)$ the free product. Let $a \in \mathcal{M}_1$, $b, c \in \mathcal{M}_2$. Assume that $a, b, c$ belong to the domain of $\sigma_{i/2}^\varphi$. Then, for every $x \in \mathcal{M}$,
\begin{equation*}
\|x - \varphi(x)1\|_\varphi \leq  \mathcal{E}(a, b, c) \max \left\{\|[x, a]\|_\varphi, \|[x, b]\|_\varphi, \|[x, c]\|_\varphi \right\} + \mathcal{F}(a, b, c)\|x\|_\varphi 
\end{equation*}
where
\begin{eqnarray*}
\mathcal{E}(a, b, c) & = & 6 \|a\|^3 + 4 \|b\|^3 + 4 \|c\|^3 \\
\mathcal{F}(a, b, c) & = & 3 \mathcal{C}(a) + 2 \mathcal{C}(b) + 2 \mathcal{C}(c) + 12 |\varphi(c b^*)| \|c b^*\| \\
\mathcal{C}(a) & = & 2 \|a\|^3 \|\sigma^\varphi_{i/2}(a) - a\| + 2 \|a\|^2 \|a^*a - 1\| \\
& & + 3(1 + \|a\|^2) \|aa^* - 1\| + 6 |\varphi(a)| \|a\|.
\end{eqnarray*}
\end{lem}

\subsection{Proof of the Theorem}

Let $\mathcal{M} := \Gamma(H_\R, U_t)''$ be a free Araki-Woods factor of type ${\rm III_1}$ and denote by $\varphi$ the free quasi-free state. We recall that such a factor can always be written as the free product of three free Araki-Woods von Neumann algebras (see the proof of Theorem 2.7 in \cite{vaes2004}):
\begin{equation*}
(\mathcal{M}, \varphi) \cong (\mathcal{M}_1, \varphi_1) \ast (\mathcal{M}_2, \varphi_2) \ast (\mathcal{M}_3, \varphi_3).
\end{equation*}
Notice that $\sigma_t^\varphi = \sigma_t^{\varphi_1} \ast \sigma_t^{\varphi_2} \ast \sigma_t^{\varphi_3}, \forall t \in \R$.

Thanks to Lemma $\ref{asymptotic}$, we may choose three sequences of unitaries $(u_n^{j})$, for $j \in \{1, 2, 3\}$, such that $u_n^j \in \mathcal{U}(\mathcal{M}_j)$ is analytic w.r.t. $(\sigma_t^{\varphi_j})$, and satisfies conditions $(1)-(3)$ of Lemma $\ref{asymptotic}$, for all $j \in \{1, 2, 3\}$. The way the sequence of unitaries $(u_n^j)$ is constructed in Lemma \ref{asymptotic} (see Lemma 4.3 in \cite{vaes2004}) shows that conditions $(1)-(3)$ are satisfied for the state $\varphi$, i.e. the sequence of unitaries $(u_n^j)$ in $\mathcal{M}_j$ satisfies for every $j \in \{1, 2, 3\}$, 
\begin{enumerate}
\item $\|\sigma_z^\varphi(u_n^j) - u_n^j\| \to 0$, uniformly on compact sets of $\C$;
\item $\varphi(u_n^j) \to 0$;
\item $\| [u_n^j, \varphi] \| \to 0$.
\end{enumerate}
Moreover, by freeness, $\varphi(u_n^3(u_n^2)^*) = \varphi(u_n^3)\overline{\varphi(u_n^2)} \to 0$. 

Assume $a \in \AB(\varphi)$. Fix $\varepsilon > 0$. Since $(u_n^j) \in \AC(\varphi)$, it follows that $[a, u_n^j] \to 0$ ultrastrongly for any $j \in \{1, 2, 3\}$, and thus we may choose $n \in \N$ large enough such that 
\begin{eqnarray*}
\|[a, u_n^j]\|_\varphi & \leq & \varepsilon/28, \forall j \in \{1, 2, 3\}, \\
\mathcal{F}(u_n^1, u_n^2, u_n^3) \|a\|_{\varphi} & \leq & \varepsilon/2.
\end{eqnarray*}
Thus, we get thanks to Lemma \ref{14eps}, $\|a - \varphi(a)1\|_{\varphi} \leq \varepsilon/2 + \varepsilon/2 = \varepsilon$. Since $\varepsilon > 0$ is arbitrary, $a = \varphi(a)1$. Thus $\AB(\varphi) = \C$ and we are done.

\subsection{Final Remark} Set $\mathcal{M} := \Gamma(L^2(\R, \R), \lambda_t)''$ the free Araki-Woods factor associated with the left regular representation $(\lambda_t)$ of $\R$ on the real Hilbert space $L^2(\R, \R)$. Shlyakhtenko showed in \cite{shlya98} that the continuous core of $\mathcal{M}$ is isomorphic to $L(\F_\infty) \otimes \mathbf{B}(\ell^2)$ and the dual action is precisely the one constructed by R\u{a}dulescu in \cite{radulescu1991}. As noticed in \cite{shlya2003}, for any f.n. state $\varphi$ on $\mathcal{M}$, the centralizer $\mathcal{M}^\varphi$ is amenable. Indeed, first we have
\begin{equation}\label{core}
\mathcal{M} \rtimes_{\sigma^\varphi} \R \simeq L(\F_\infty) \otimes \mathbf{B}(\ell^2).
\end{equation}
Choose on the left hand-side of $(\ref{core})$ a non-zero projection $p \in L(\R)$ such that $\Tr(p) < +\infty$. We know that $p(\mathcal{M} \rtimes_{\sigma^\varphi} \R)p \simeq L(\F_\infty)$ is {\it solid} by Ozawa's result (\cite{ozawa2003}). Since $L(\R)p$ is diffuse in $p(\mathcal{M} \rtimes_{\sigma^\varphi} \R)p$, its relative commutant must be amenable. In particular $\mathcal{M}^\varphi \otimes L(\R)p$ is amenable. Thus, $\mathcal{M}^\varphi$ is amenable. Consequently, we obtain 

\begin{cor}
Let $\mathcal{M} := \Gamma(L^2(\R, \R), \lambda_t)''$. Then there exists a f.n. state $\psi$ on $\mathcal{M}$ such that $(\mathcal{M}^\psi)' \cap \mathcal{M} = \C$. Moreover, $\mathcal{M}^\psi$ is isomorphic to the unique hyperfinite ${\rm II_1}$ factor.
\end{cor}

\bibliographystyle{plain}

\end{document}